\documentclass[10pt]{amsart}
\usepackage{amsmath,amsthm,amscd,amssymb}
\usepackage{graphicx}
\usepackage[dvipsnames,usenames]{color}
\usepackage{subfigure}

\usepackage[colorlinks=true, urlcolor=NavyBlue, linkcolor=NavyBlue, citecolor=NavyBlue]{hyperref}

\newtheorem{theorem}{Theorem}[section]

\newtheorem{corollary}[theorem]{Corollary}

\newtheorem{lemma}[theorem]{Lemma}

\newtheorem{proposition}[theorem]{Proposition}

\numberwithin{equation}{section}
\numberwithin{figure}{section}

\newtheorem*{thm:mainappliction}{Theorem~\ref{}}

\theoremstyle{definition}

\newtheorem{question}[theorem]{Question}

\newcommand{\Z}{\mathbb{Z}}

\newcommand{\Q}{\mathbb{Q}}

\author{Tim D. Cochran$^\dag$}
\address{Department of Mathematics\\Rice University\\P.O. Box 1892\\Houston, TX 77251}
\email{cochran@math.rice.edu}
\urladdr{http://math.rice.edu/~cochran}
\author{Bridget D. Franklin$^{\dag\dag}$}
\address{Department of Mathematics\\Rice University\\P.O. Box 1892\\Houston, TX 77251}
\email{bridget.franklin@rice.edu}
\urladdr{http://math.rice.edu/~bf2}
\author{Peter D. Horn$^{\dag\dag\dag}$}
\address{Department of Mathematics\\Columbia University\\2990 Broadway\\New York, NY 10025}
\email{pdhorn@math.columbia.edu}
\urladdr{http://math.columbia.edu/~pdhorn}

\thanks{$^{\dag}$Partially supported by the National Science Foundation  DMS-1006908}
\thanks{ $^{\dag\dag}$Partially supported by Nettie S. Autry Fellowship}
\thanks{$^{\dag\dag\dag}$Partially supported by NSF Postdoctoral Fellowship DMS-0902786}

\title{Rational Knot Concordance and Homology Cobordism}

\begin{document}

\begin{abstract} The following is a long-standing open question: `` If the zero-framed surgeries on two knots in $S^3$ are $\mathbb{Z}$-homology cobordant, are the knots themselves concordant?''  We show that an obvious rational version of this question has a negative answer. Namely, we give examples of knots whose zero-framed surgeries are $\mathbb{Q}$-homology cobordant $3$-manifolds,  wherein the knots are not  $\mathbb{Q}$-concordant (that is not concordant in any rational homology $S^3\times [0,1]$). Specifically, we prove that, for any positive integer $p$ and any knot $K$, the zero framed surgery on $K$ is $\Z\left[\frac{1}{p}\right]$-homology cobordant to the zero framed surgery on its $(p,\pm 1)$ cables. Then we observe that most knots are not $\mathbb{Q}$-concordant to their $(p,\pm 1)$ cables.
\end{abstract}
\maketitle
\section{Introduction}\label{sec:intro}

A \textbf{knot} $K$ is the image of a tame embedding of a circle into $S^3$. Two knots, $K_0\hookrightarrow S^3\times \{0\}$ and $K_1\hookrightarrow S^3\times \{1\}$, are \textbf{CAT concordant} (CAT= SMOOTH or TOPOLOGICAL LOCALLY FLAT) if there exists a proper CAT embedding of an annulus into $S^3\times [0,1]$ that restricts to the knots on $S^3\times \{0,1\}$. A \textbf{CAT slice knot} is one that is the boundary of a CAT embedding of a $2$-disk in $B^4$. Recall that to any knot, $K$, is associated a closed $3$-manifold, $M_K$, called the \emph{zero framed surgery on $K$}, which is obtained from $S^3$ by removing a solid torus neighborhood of $K$ and replacing it with another solid torus in such a way that the longitude of $K$ bounds the meridian of the new solid torus. The question of whether two knots, $K_0$ and $K_1$, are CAT concordant is closely related to the question of whether or not the associated $3$-manifolds, $M_{K_0}$ and $M_{K_1}$, are CAT homology cobordant. For if the knots are concordant via a CAT embedded annulus $A$ then, by Alexander duality, the exterior of $A$ is a homology cobordism relative boundary between the exteriors of the knots. It follows, by adjoining $S^1\times D^2\times [0,1]$, that the zero framed surgeries are homology cobordant (i.e. one does ``zero framed surgery on the annulus''). The converse is a long-standing open question:

\begin{question}\label{question:Kirby} If $M_{K_0}$ and $M_{K_1}$ are homology cobordant (if necessary also assume that $\pi_1$ of the cobordism is normally generated by either meridian) then are $K_0$ and $K_1$ concordant?  (See Kirby's Problem List ~\cite[Problem $1.19$ and remarks]{Kproblemlist}).
\end{question}

Note that Question~\ref{question:Kirby} has a negative answer if knots are assumed to be oriented.  Livingston has provided examples of oriented knots which are not concordant to their reverse \cite{Li7}, though the zero framed surgeries are actually homeomorphic.  

Evidence for a positive answer to Question~\ref{question:Kirby} in the more general case is provided by the following well-known result.

\begin{proposition}\label{prop:basichomcobordism}  Suppose $K$ is a knot and $U$ is the trivial knot. Then the following are equivalent:
\begin{itemize}
\item [a.] $M_K$ is CAT homology cobordant to $M_U$ via a cobordism $V$ whose $\pi_1$ is normally generated by a meridian of $K$.
\item [b.] $M_K=\partial W$, where $W$ is a CAT manifold  that is a homology circle and whose $\pi_1$ is normally generated by the meridian.
\item [c.] $K$ bounds a CAT embedded $2$-disk in a CAT manifold $\mathcal{B}$ that is homeomorphic to $B^4$.
\item [d.] $K$ is CAT concordant to $U$ in a CAT $4$-manifold that is homeomorphic to $S^3\times [0,1]$.
\end{itemize}
Moreover $d \Rightarrow a$ even if $U$ is not assumed to be the trivial knot.
\end{proposition}

In particular the case $a \Rightarrow d$ of Proposition~\ref{prop:basichomcobordism} in the TOP category yields:
\begin{corollary}\label{cor:truTOP} Question~\ref{question:Kirby} has a positive answer in the topological category if one of the knots is the trivial knot. The same is true in the smooth category if $B^4$ has a unique smooth structure up to diffeomorphism.
\end{corollary}

\begin{proof}[Proof of Proposition~\ref{prop:basichomcobordism}] The implication $d\Rightarrow a$ has already been explained above (and is true for any two knots). Note that $M_U\cong S^1\times S^2$, which is the boundary of the homology circle $S^1\times B^3$. Adjoining the latter to one end of the homology bordism $V$, provided by $a.$, yields $W$ which shows $a\Rightarrow b$. Given $W$ from $b.$, add a zero-framed $2$-handle to $\partial W$ along the meridian of $K$. The resulting manifold $\mathcal{B}$ is a simply-connected homology $4$-ball whose boundary is $S^3$. By Freedman's theorem, $\mathcal{B}$ is homeomorphic to $B^4$. The core of the attached $2$-handle is a flat disk whose boundary is a copy of $K\hookrightarrow S^3$. Thus $b\Rightarrow c$. To show $c\Rightarrow d$, merely remove a small $4$-ball from $\mathcal{B}$ centered at a point of the $2$-disk.
\end{proof}

Remarkably, Proposition~\ref{prop:basichomcobordism} has an analogue for \emph{rational} homology cobordism. In particular it is related to notions of \emph{rational concordance} which have also been  previously studied \cite{CO1,ChaKo1,ChaKo2,Gi6, Ka4} and finally treated systematically by Cha in ~\cite{Cha2}. Before stating this analogue, we review some terms.

Suppose that $R\subset\Q$ is a non-zero subring. Recall that a space $X$ \emph{is an $R$-homology $Y$} means that $H_*(X;R)\cong H_*(Y;R)$.  Knots $K_0$ and $K_1$ in $S^3$ are said to be \emph{CAT $R$-concordant} if there exists a compact, oriented CAT $4$-manifold $W$, that is an $R$-homology $S^3\times [0,1]$, whose boundary is $S^3\times \{0\}\sqcup -(S^3\times \{1\})$, and in which there exists a properly CAT embedded annulus $A$ which restricts on its boundary to the given knots. We then say that $K_0$ is \emph{CAT $R$-concordant} to $K_1$ \emph{in $W$}. We say that $K$ is \emph{CAT $R$-slice} if it is CAT $R$-concordant to $U$ or, equivalently, if it bounds a CAT embedded $2$-disk in an $R$-homology $4$-ball whose boundary is $S^3$. The latter notion agrees with ~\cite{ChaKo1,Cha2} but is what Kawauchi calls weakly $\Q$-slice \cite{Ka4}.

\begin{proposition}\label{prop:Rbasichomcobordism}  Suppose $K$ is a knot in $S^3$ and $R\subset\Q$ is a non-zero subring. Then the following are equivalent:
\begin{itemize}
\item [a.] $M_K$ is CAT $R$-homology cobordant to $M_U$.
\item [b.] $M_K=\partial W$, where $W$ is a CAT manifold  that is an $R$-homology circle.
\item [c.] $K$ is CAT $R$-slice.
\item [d.] $K$ is CAT $R$-concordant to $U$.
\end{itemize}
Moreover $d \Rightarrow a$ for any two knots.
\end{proposition}

The proof of Proposition~\ref{prop:Rbasichomcobordism} is essentially identical to that of Proposition~\ref{prop:basichomcobordism}.

Proposition~\ref{prop:Rbasichomcobordism} suggests that invariants that obstruct knots from being rationally concordant might be dependent only on the rational homology cobordism class of the zero framed surgery. Further evidence for this is  by the following observation,  which follows immediately from $d \Rightarrow c$ of Proposition~\ref{prop:Rbasichomcobordism} and Ozsvath-Szabo~\cite[Theorem 1.1]{OzSz2}.

\begin{corollary}\label{cor:tauratU} For any $R$, if $M_K$ is smoothly $R$-homology cobordant to $M_U$, then $\tau(K)=0$.
\end{corollary}

This suggests the following ``rational version'' of Question~\ref{question:Kirby}, which is the question of whether or not  $a \Rightarrow d$ of Proposition~\ref{prop:Rbasichomcobordism} holds for any two knots (for $R=\Q$).

\begin{question}\label{question:KirbyQ} If $M_{K_0}$ and $M_{K_1}$ are CAT $\Q$-homology cobordant,  then are $K_0$ and $K_1$ CAT $\Q$-concordant?
\end{question}

\noindent Of course by Proposition~\ref{prop:Rbasichomcobordism} this question has a positive answer if one of the knots is the unknot.

In this paper, we do not resolve Question~\ref{question:Kirby}, but show that Question~\ref{question:KirbyQ} has a \emph{negative} answer. To accomplish this we first prove in Section~\ref{sec:mainresult} :

\newtheorem*{prop:mainresult}{Proposition~\ref{prop:mainresult}}
\begin{prop:mainresult} For any knot $K$ and any positive integer $p$, zero framed surgery on $K$ is smoothly $\Z\left[\frac{1}{p}\right]$-homology cobordant to zero framed surgery on the $(p,1)$-cable of $K$ (and also to the zero framed surgery on the $(p,-1)$-cable of $K$).
\end{prop:mainresult}

\noindent Then in Section~\ref{sec:Qinvts} we observe that there are elementary classical invariants that obstruct a knots being topologically rationally concordant to its $(p,1)$-cable. Even among topologically slice knots, the $\tau$ invariant can be used to obstruct $K$  being smoothly rationally concordant to its $(p,1)$-cable, and we do so, using work of M. Hedden. 

Thus in summary we show:
\newtheorem*{thm:mainapplication}{Theorem~\ref{thm:mainapplication}}
\begin{thm:mainapplication} The answer to Question~\ref{question:KirbyQ} is ``No,'' in both the smooth and topological category. In the smooth category there exist counterexamples that are topologically slice.
\end{thm:mainapplication}

Thus the answer to Question~\ref{question:KirbyQ} is decidedly ``No.''  In the final section, we formulate a refined version of Question~\ref{question:KirbyQ},  and this refined version, like Question~\ref{question:Kirby}, remains open.

We remark in passing that the analogues of Proposition~\ref{prop:basichomcobordism} and Corollary~\ref{cor:truTOP} ~\cite[Theorem 5.2, p.19]{Hi} hold for links and have been the basis of recent attempts to resolve the smooth $4$-dimensional Poincar\'{e} Conjecture ~\cite{FGMW}.  However, Question~\ref{question:Kirby} has a negative answer for links. Cochran-Orr gave examples in ~\cite[Fig.1]{CO1a}~\cite{CO1} of links whose zero surgeries are \emph{diffeomorphic} but which are not even TOP concordant. There are even such examples with distinct Milnor's $\overline{\mu}$-invariants. (Consider the link, $L_m$,  whose first component is an $m$-twist knot and whose second component is a trivial circle linking the $m$-twisted band. The zero framed surgery on $L_m$ is then independent of $m$ but the concordance type of the first component is dependent of $m$!)

\section{Rational homology cobordism and cable knots}\label{sec:mainresult}

\begin{proposition}\label{prop:mainresult} For any knot $K$ and any positive integer $p$, zero framed surgery on $K$ is smoothly $\Z\left[\frac{1}{p}\right]$-homology cobordant to zero framed surgery on the $(p,1)$-cable of $K$ (and also to the zero framed surgery on the $(p,-1)$-cable of $K$).
\end{proposition}

We present two proofs for the case of the $(p,1)$-cable. The second proof also clearly covers the $(p,-1)$ case.

\begin{proof}[Proof of Proposition~\ref{prop:mainresult} for $p=2$] 

	In this proof we freely use the basic techniques of the calculus of framed links, commonly called Kirby calculus, that can be used to encode handlebody descriptions of $3$- and $4$-dimensional manifolds ~\cite{K}. Let $M_K$ denote the 3-manifold obtained by $0$-framed Dehn filling on $K$. Then $H_1(M_K)\cong \Z$, generated by a meridian of $K$. Consider $M_K\times [0,1]$. Attach to $M_K\times \{1\}$ a 1- and a 2-handle according to Figure~\ref{addhandles}, and call the resulting 4-manifold $W$.  The 1-handle adds a copy of $\Z$ to $H_1$, and the 2-handle equates two times a generator of this $\Z$ factor with the meridian of $K$, denoted $\mu_K$.  Thus $W$ has the same $\Z\left[\frac{1}{2}\right]$-homology as $M_K$.
	
	\begin{figure}[hbt!]
		\begin{center}
			\begin{picture}(269,63)(0,0)
				\put(6, 28){$K$}
				\put(40, 5){$0$}
				\put(85, 28){$\leadsto$}
				\put(113, 28){$K$}
				\put(147, 5){$0$}
				\put(193, 5){$0$}
				\includegraphics{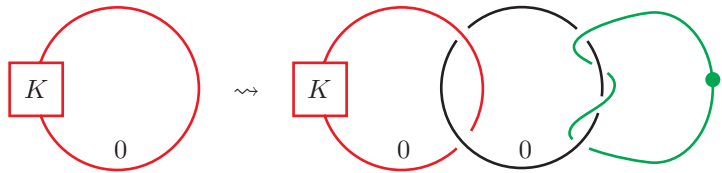}
			\end{picture}
			\caption{Add a 1- and 2-handle.}
			\label{addhandles}
		\end{center}
	\end{figure}
	
\noindent	The `bottom' boundary of $W$ is $-M_K$.  The `top' boundary component is $M_{K(2,1)}$.  To see this, note that the `top' boundary component is diffeomorphic to the 3-manifold in Figure~\ref{timsmoves}(a).  To reach the diagram in Figure~\ref{timsmoves}(b), slide the green (right-most) handle over the red (left-most) handle twice, as shown.
	
	\begin{figure}[h!]
	    \begin{center}
	    	\begin{picture}(350,69)(0,0)
	    		\subfigure[Step one]
			    {
			        \put(6, 28){$K$}
					\put(40,5){$0$}
					\put(84,5){$0$}
					\put(135, 6){$0$}
					\includegraphics{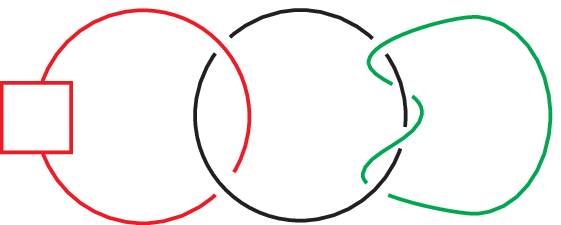}
			        \label{move1}
			    }
				\hspace{1cm}
			    \subfigure[Step two]
			    {
					\put(6, 36){$K$}
					\put(41, 11){$0$}
					\put(85, 11){$0$}
					\put(136, 12){$0$}
			        \includegraphics{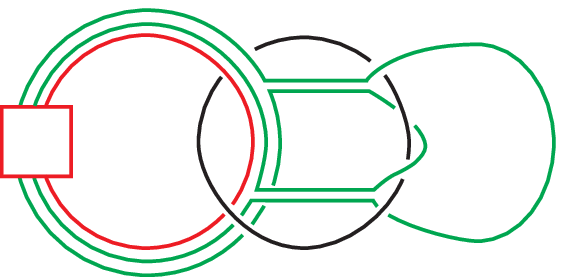}
			        \label{move2}
			    }
	    	\end{picture}
		    \caption{Bottom boundary.}
		    \label{timsmoves}
	    \end{center}
	\end{figure}
	
\noindent	Now isotope the green (right-most) handle in Figure~\ref{timsmoves}(b) to arrive at Figure~\ref{timsmoves2}(a).  Isotope further to arrive at Figure~\ref{timsmoves2}(b), where the black (small) curve is isotopic to a meridian of $K$.  We may use this curve as a helper circle to unlink the red (innermost) curve from the green and then undo crossings of the red curve, arriving at Figure~\ref{timsmoves2}(c).  The 3-manifold described in Figure~\ref{timsmoves2}(c) is diffeomorphic to $M_{K(2,1)}$.  Thus, $\partial W$ is the disjoint union of $M_{K(2,1)}$ and $-M_K$.
	
	\begin{figure}[h!]
	    \begin{center}
	    	\begin{picture}(450,69)(0,0)
	    		\subfigure[Step three]
			    {
			        \put(6, 35){$K$}
					\put(40,12){$0$}
					\put(84,12){$0$}
					\put(135, 13){$0$}
					\includegraphics{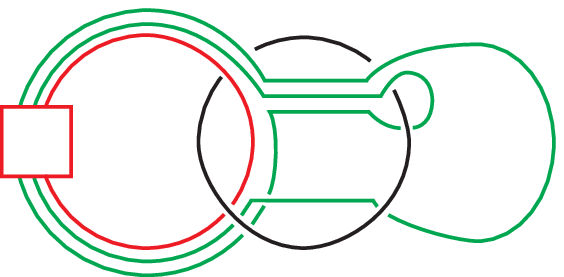}
			        \label{move3}
			    }
				\hspace{.5cm}
			    \subfigure[Step four]
			    {
					\put(6, 38){$K$}
					\put(41, 13){$0$}
					\put(80, 13){$0$}
					\put(55, 28){$0$}
			        \includegraphics{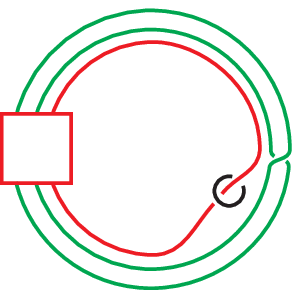}
			        \label{move4}
			    }
				\hspace{.5cm}
				\subfigure[Step five]
			    {
					\put(6, 38){$K$}
					\put(82, 13){$0$}
					\put(93, 28){$0$}
					\put(128, 28){$0$}
			        \includegraphics{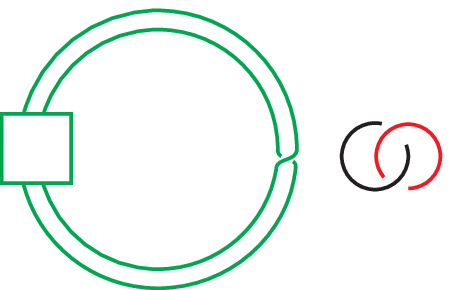}
			        \label{move5}
			    }
	    	\end{picture}
		    \caption{Isotoping handles.}
		    \label{timsmoves2}
	    \end{center}
	\end{figure}		
\end{proof}

The proof for $p>2$ is similar but requires the attaching circle of the black 2-handle we added to wrap $p$ times around the 1-handle and then requires $p$ handleslides.  Instead of drawing these pictures, we provide a different proof for the general case.

\begin{proof}[Another proof of Proposition~\ref{prop:mainresult}]
 	Let $U=U(p,1)$ denote the $(p,1)$-torus knot and let $\alpha$ denote a push-off of the meridian of the solid torus on which $U(p,1)$ lies (so that $\ell k(\alpha,U)=p)$. Thus in homology $[\alpha]=p\,[\mu_{U}]$. Recall that the knot $K(p,1)$ is a satellite knot with $U(p,1)$ as pattern and companion $K$. As such the exterior of $K(p,1)$ decomposes into two pieces, the exterior of $K$ and the exterior of $U(p,1)\cup \alpha$, glued along the (toral) boundary of $\nu(\alpha)$. Thus the $0$-framed Dehn filling, $M_{K(p,1)}$, decomposes (along the same torus) into $S^3\setminus K$ and $M_U\setminus \nu(\alpha)$. Additionally recall that $M_K$ decomposes into $S^3\setminus \nu(K)$ and a solid torus.

With these facts in mind, we now construct  a $4$-manifold $E$ whose boundary is the disjoint union $M_K \sqcup M_{U} \sqcup -M_{K(p,1)}$. Begin with the disjoint union of $M_K  \times [0, 1]$ and $M_{U} \times [0,1]$. Identify the solid torus $\nu(\alpha) \times \{1\}$ in $M_{U} \times \{1\}$ with the solid torus $(M_K \setminus (S^3 - \text{int}\ \nu(K)) \times \{1\}$ in  $M_K  \times [0, 1]$. Do this in such a way that (a parallel push-off of) $\alpha$ is identified to a meridian of $K$ (which is a longitude of the latter solid torus). Then the third boundary component of $E$ is
$$
\left(M_U\setminus \nu(\alpha)\right)\cup S^3\setminus K\equiv M_{K(p,1)},
$$
as claimed.
Furthermore, note that under the inclusion maps on homology,
\begin{equation}\label{eq:h1facts}
[\mu_K\times \{0\}]\sim[\mu_K\times \{1\}]\equiv [\alpha\times \{1\}]\sim [\alpha\times \{0\}]\sim p\,[\mu_{U}],
\end{equation}
where for the last we use our remark in the first paragraph above. We may analyze the homology of $E$ by the Mayer-Vietoris sequence below
$$
H_2(\nu(\alpha))\to H_2(M_U)\oplus H_2(M_K)\to H_2(E)\to H_1(\nu(\alpha))\overset{\psi}{\rightarrow} H_1(M_U)\oplus H_1(M_K)\to H_1(E)\to 0
$$
By ~\eqref{eq:h1facts}, $\psi$ is injective with infinite cyclic cokernel. Since $H_2(\nu(\alpha))=0$, this yields the following elementary lemma which may be compared with \cite[Lemma 2.5]{CHL3}.

\begin{lemma}\label{lem:Efacts}
The inclusion maps induce the following
\begin{enumerate}	
	\item an isomorphism $H_1\left(M_{U};\mathbb{Z}\right)\cong H_1\left(E;\mathbb{Z}\right)$;
\item an isomorphism $H_1\left(M_{K};\mathbb{Z}\left[\frac{1}{p}\right]\right)\rightarrow H_1\left(E;\mathbb{Z}\left[\frac{1}{p}\right]\right)$;
	\item an isomorphism $H_2(E;\mathbb{Z}) \cong H_2\left(M_{U};\mathbb{Z}\right)\oplus H_2(M_{K};\mathbb{Z})$
		\end{enumerate}
	\end{lemma}

	Finally, let $S = B^4\setminus \Delta$ where $\Delta\subset B^4$ is a slice disk for the unknot $U(p,1)$. Then $S$ is an integral homology circle whose first homology is generated by $\mu_U$ and whose boundary is $M_U$. Let $W$ be the space obtained by attaching $S$ to $E$ along $M_{U}$.  $W$ is a 4-manifold with boundary $M_K \sqcup -M_{K(p,1)}$.

	Another Mayer-Vietoris argument gives
$$
H_2(W;\mathbb{Z}) \cong \frac{H_2(E)}{H_2(M_U)}\cong H_2(M_K;\mathbb{Z}),
$$
using ($3$) of Lemma~\ref{lem:Efacts}; and
$$
H_1(W;\mathbb{Z}) \cong H_1(E;\mathbb{Z}).
$$
Combining these facts we conclude that
$$
0=H_2\left(W,M_K;\mathbb{Z}\left[1/p\right]\right)\cong H_1\left(W,M_K;\mathbb{Z}\left[1/p\right]\right)\cong H_1\left(W,M_{K(p,1)};\mathbb{Z}\left[1/p\right]\right).
$$
It the follows from Poincare duality that $W$ is a $\mathbb{Z}[\frac{1}{p}]$-homology cobordism between $M_K$ and $M_{K(p,1)}$.
\end{proof}

\section{Invariants of rational concordance}\label{sec:Qinvts}

In this section we will observe that many classical concordance invariants obstruct knots being $\Q$-concordant. This question has been considered in even greater generality in ~\cite{Cha2}. From this we deduce that only rarely is a knot $K$  $\Q$-concordant to its cable $K(p,1)$. We then observe that the $\tau$-invariant of Oszvath-Szabo and Rasmussen can be used to obstruct smooth $\Q$-concordance (even between topologically slice knots). This is then used, in conjunction with known computations of $\tau$, to give examples of topologically slice knots $K$ which are not $\Q$-concordant to \emph{any} of of their cables $K(p,1)$.

Before beginning we should point out that there \emph{do} exist non-slice knots $K$ for which $K$ is smoothly $\Q$-concordant to $K(p,1)$ for \emph{every} non-zero $p$. For suppose that $K$ is a smoothly $\Q$-slice knot that is not a smoothly slice knot. Examples of this are provided by the figure-eight knot (\cite[p.63]{Cha2}\cite[Lemma 2.2]{CHL6}), or more generally any (non-slice) strongly negative-amphicheiral knot ~\cite{Ka4}. Such a knot $K$ is $\Q$-concordant to the unknot $U$. Thus $K(p,1)$ is $\Q$-concordant to the unknot $U(p,1)$ for any $p$. Thus $K$ is $\Q$-concordant to $K(p,1)$.

\begin{proposition}\label{prop:coolexs} Suppose $K$ is a strongly negative-amphicheiral knot. Then for every non-zero $p$, $K$ is $\Q$-concordant to $K(p,1)$.
\end{proposition}

However, this is rare as the following subsection shows.

\begin{subsection}{Topological Invariants of rational concordance}

If $K$ is a knot in $S^3$ and $V$ is a Seifert matrix for $K$ then recall the Levine-Tristram $\omega$-signature of $K$, $\sigma_K(\omega)$, for any $|\omega|=1$ is the signature of
$$
(1-\omega)V+(1-\overline{\omega})V^T.
$$
The only explicit reference we found for the following theorem is ~\cite[Thm. 1.1]{ChaKo1} (see also ~\cite[Theorem 1.7]{CO1}\cite{CO1a}).

\begin{proposition}\cite[Thm. 1.1]{ChaKo1}\label{prop:sigsinvt} If $K_0$ is CAT $\Q$-concordant to $K_1$ then for any prime $p$ and $\omega=e^{2\pi ik/p}$,
$$
\sigma_{K_0}(\omega)=\sigma_{K_1}(\omega).
$$
\end{proposition}

Indeed Cha shows that the ordinary algebraic knot concordance group embeds into an algebraic $\Z_{(2)}$-concordance group ~\cite[Section 2.2, 4.4]{Cha2} ( see also ~\cite[Theorem 1.7]{CO1}\cite{CO1a}).

\noindent Note that, since such values of $\omega$ are dense in the circle and since the Levine signature functions of $K_0$ and $K_1$ are constant except at roots of the Alexander polynomial, Proposition~\ref{prop:sigsinvt} implies that these functions agree for $K_0$ and $K_1$, except possibly at roots of their respective Alexander polynomials.

\begin{corollary}\label{cor:finiteorder} If $K$ is CAT $\Q$-concordant to $K(p,1)$ for some $p>1$, then $K$ is of finite order in the algebraic knot concordance group.
\end{corollary}
\begin{proof} [Proof of Corollary~\ref{cor:finiteorder}] Assume that $K$ is TOP $\Q$-concordant to $K(p,1)$. Then their signatures agree (except possibly at roots of the Alexander polynomials). Suppose \emph{some} Levine-Tristram-signature of $K$ were non-zero. Since the $\omega=0$ signature always vanishes and since the Levine-Tristram signature function is locally constant, possibly jumping only at roots of the Alexander polynomial, we can  choose $\omega$ so that $\omega^p$ is (or is close to) the ``first'' value on the unit circle (smallest argument) for which $K$ has non-zero signature (and avoiding roots of the Alexander polynomial). That is we can choose $\omega$ so that
$$
\sigma_{K}(\omega)=0, ~~\text{and}~~ \sigma_{K}(\omega^p)\neq 0.
$$
But it is known by ~\cite{Kea2} that,
$$
\sigma_{K(p,1)}(\omega)=\sigma_K(\omega^p).
$$
Combining this with Proposition~\ref{prop:sigsinvt}, we see that
$$
\sigma_K(\omega^p)=\sigma_K(\omega).
$$
This is false for our particular choice of $\omega$ above. Hence the signature function of $K$ vanishes (excluding roots of the Alexander polynomial). It is known that this is equivalent to $K$ being of finite order in the algebraic knot concordance group ~\cite{L5}.
\end{proof}

\begin{corollary}\label{cor:trefoil1} The right-handed trefoil knot, $T$, is not TOP $\Q$-concordant to $T(p,1)$ for any $p>1$.
\end{corollary}

There are other related papers that discuss the question of whether a knot is $\Z$-concordant to its $(p,1)$-cable ~\cite{Ka5}\cite{Lith1}\cite{CRL}.

The following criteria can be applied even when the knot signatures fail. The first does not seem to appear in the literature although it does follow, for example, from combining results of the much more general \cite{Cha2}. The second is implicit in \cite{Cha2}. We sketch a proof in order to make a pedagogical point about rational concordance.

\begin{proposition}\label{prop:polysinvt} If $K_0$ is CAT $\Q$-concordant to $K_1$ then for some positive integer $k$ and for some integral polynomial  $f$,
$$
\delta_{0}\left(t^k\right)\delta_{1}\left(t^k\right)\doteq f(t)f\left(t^{-1}\right),
$$
where $\delta_{i}(t)$ is the Alexander polynomial of $K_i$.
\end{proposition}

Most generally, let $\mathcal{B}\ell^K(t)$ denote the nonsingular Blanchfield linking form defined on the rational Alexander module of $K$, $\mathcal{A}^\Q(K)\equiv H_1(M_K;\Q[t,t^{-1}])$. Then,

\begin{proposition}\label{prop:Blanchinvt} If $K_0$ is CAT $\Q$-concordant to $K_1$ then for some positive integer $k$,
$$
\mathcal{B}\ell^{K_0}(t^k)\sim \mathcal{B}\ell^{K_1}(t^k),
$$
where these denote the induced forms on the module
$$
\mathcal{A}^\Q(K_i)\otimes_{\Q[t,t^{-1}]}\Q[t,t^{-1}]
$$
where here the right-hand $\Q[t,t^{-1}]$ is a module over itself via the map $t\to t^k$; and $\sim$ denotes equality in the Witt group of such forms (see \cite{Cha2}\cite{Hi}).
\end{proposition}

\begin{proof} Suppose $K_0$ is CAT $\Q$-concordant to $K_1$ via an annulus, $A$, embedded in a $\Q$-homology $S^3\times [0,1]$, $W$. Let $E_0$, $E_1$ and $E_A$ denote the exteriors of $K_0$, $K_1$ and $A$ respectively. Then
$$
\frac{H_1(E_*;\Z)}{\mathrm{torsion}}\cong \Z
$$
for $*=0,1,A$. The \textbf{complexity of the concordance}  is the positive integer, $k$, for which the image of the meridian $\mu_i$, for $i=0,1$,   under the inclusion-induced map $j_i$
$$
\frac{H_1(E_0;\Z)}{\mathrm{torsion}}\overset{j_0}{\longrightarrow}\frac{H_1(E_A;\Z)}{\mathrm{torsion}}
\overset{j_1}{\longleftarrow}\frac{H_1(E_1;\Z)}{\mathrm{torsion}},
$$
is $\pm k$ times a generator. This was defined in \cite{CO1,CO1a} but see also \cite{Cha2, ChaKo1} and was called the \emph{multiplicity} in \cite[page 463]{COT}. There is a unique epimorphism
$$
\phi:\pi_1(E_A)\to \Z.
$$
This defines a coefficient system on $E_A$ and also on $E_i$ for $i=0,1$  by setting $\phi_i=\phi\circ j_i$. Then it is well-known that the Alexander modules using these induced coefficient systems are not the ordinary Alexander modules but rather,
$$
H_1(E_i;\Z[t,t^{-1}])\cong \mathcal{A}(K_i)\otimes_{\Z[t,t^{-1}]}\Z[t,t^{-1}]
$$
where the right-hand $\Z[t,t^{-1}]$ is a module over itself via the map $t\to t^k$. The order of such a module is well-known to be $\delta_{i}(t^k)$ where $\delta_{i}(t)$ is the order of $\mathcal{A}(K_i)$. (This ``tensored up''  module is the same as the Alexander module of the $(k,1)$-cable of $K_i$). The coefficient system $\phi$ also induces Blanchfield linking forms on these modules and these differ from the ordinary Blanchfield form in the analogous manner.

If $A$ were an actual concordance then we have the classical result that the kernel, $P$, of the map
$$
\mathcal{M}\equiv H_1(E_0;\Q[t,t^{-1}])\oplus H_1(E_1;\Q[t,t^{-1}])\to H_1(E_A;\Q[t,t^{-1}])
$$
is self-annihilating with respect to the ordinary Blanchfield forms. It would then follow (by definition) that the Blanchfield forms are equivalent in the Witt group. It also would then follow that $P$ is isomorphic to the dual of $\mathcal{M}/P$, quickly yielding the classical result 
$$
\delta_{K_0}\left(t\right)\delta_{K_1}\left(t\right)\doteq f(t)f\left(t^{-1}\right),
$$
for some polynomial $f$. In the situation that $A$ is only a $\Q$-concordance, these results are also known (see for example \cite[Theorem 4.4, Lemma 2.14]{COT}). The only difference is that the relevant modules and forms are not the ordinary Alexander modules but rather are  ``tensored up'' as above; and the orders of the relevant modules are not the actual Alexander polynomials of $K_i$, but are $\delta_i(t^k)$. The claimed results follow.
\end{proof}

\begin{corollary}\label{cor:twist} Suppose $K$ is the $3$-twist knot with a negative clasp . Then, although $K$ is of finite order in the algebraic concordance group, $K$ is not TOP $\Q$-concordant to $K(p,1)$ for any $p>1$.
\end{corollary}
\begin{proof} The Alexander polynomial of $K$ is $\delta(t)=3t-7+3t^{-1}$, whereas the Alexander polynomial of $K(p,1)$ is $\delta(t^p)$. If $K$ were TOP $\Q$-concordant to $K(p,1)$  then by Proposition~\ref{prop:polysinvt}, for some positive $k$ and integral polynomial $f(t)$,
$$
\delta\left(t^k\right)\delta\left(t^{kp}\right)=\pm t^g f(t)f\left(t^{-1}\right).
$$
But $\delta\left(t^k\right)$ (and thus ($\delta\left(t^{kp}\right)$) is irreducible for any $k$ ~\cite[Proposition 3.18]{Cha2}. This  contradicts unique factorization  if $p>1$.
\end{proof}

Casson-Gordon invariants and higher-order von-Neumann signatures should yield higher-order obstructions to $\Q$-concordance.

\end{subsection}
\begin{subsection}{Smooth Rational concordance invariants for topologically slice knots}

The Ozsvath-Szabo-Rasmussen $\tau$-invariant is an integral-valued knot invariant that is invariant under smooth concordance and additive under connected sum ~\cite{OzSz2}. It is not invariant under topological concordance and therefore may be used in cases where algebraic invariants fail. It is also known  that it is an invariant of smooth rational concordance.

\begin{proposition}\label{prop:tau1} If $K$ is smoothly $R$-concordant to $J$ then $\tau(K)=\tau(J)$.
\end{proposition}
\begin{proof} We are given that $K$ and $J$ are connected by a smooth annulus $A$ in a smooth $R$-homology $S^3\times [0,1]$, $W$. Choose an arc in $A$ from $K$ to $J$. By deleting a small neighborhood of this arc from $W$ we arrive at a smooth $R$-homology $4$-ball $\mathcal{B}$. The annulus $A$ is cut open yielding a $2$-disk whose boundary is the knot type of $K\#-J$. Thus $K\#-J$ is smoothly $R$-slice.  By ~\cite[Theorem 1.1]{OzSz2}, $\tau(K\#-J)=0$, so $\tau(K)=-\tau(-J)=\tau(J)$, the last property being also established in ~\cite{OzSz2}. 
\end{proof}

\end{subsection}

\section{The Main Result}\label{sec:mainthm}

\begin{theorem}\label{thm:mainapplication} The answer to Question~\ref{question:KirbyQ} is ``No,'' in both the smooth and topological category. In the smooth category there exist counterexamples that are topologically slice. 
\end{theorem}
\begin{proof}[Proof of Theorem~\ref{thm:mainapplication}] Let $T$ be the trefoil knot, or indeed any knot with some non-zero Levine-Tristram signature. By Proposition~\ref{prop:mainresult}, for any $p>1$, $M_T$ is (TOP and SMOOTH) $\Q$-homology cobordant to $M_{T(p,1)}$. But  by Corollary~\ref{cor:trefoil1}, or more generally by Corollary~\ref{cor:finiteorder}, $T$ is neither TOP nor SMOOTH $\Q$-concordant to $T(p,1)$. Therefore the answer to Question~\ref{question:KirbyQ} is ``No.'' Proposition~\ref{prop:polysinvt} can also be used to give examples that have finite order in the algebraic knot concordance group.

Now consider the smooth category. Let $K_0$ be the untwisted, positively-clasped Whitehead double of the right-handed trefoil knot, and let $K_1$ be the $(p,1)$-cable of $K_0$.  The Alexander polynomials of $K_0$ and $K_1$ are  equal to $1$, and so work of Freedman~\cite{FQ} implies these knots are topologically slice. The zero-framed surgeries on these knots are smoothly $\mathbb{Q}$-homology cobordant by Proposition~\ref{prop:mainresult}. If these knots were smoothly $\mathbb{Q}$-concordant then by Proposition~\ref{prop:tau1}, $\tau(K_0)=\tau(K_1)$. But this is not true. In ~\cite{He1} Hedden proves $\tau(K)=1$ whereas in ~\cite[Theorem 1.2]{He2}, he shows $\tau(K(p,1))=p\,\tau(K)=p$.  Thus if $p\neq 1$, these knots are not smoothly $\mathbb{Q}$-concordant. Therefore  there exist topologically slice knots for which the answer to Question~\ref{question:KirbyQ} is ``No'' in  the smooth category.
\end{proof}

\section{Other questions}\label{sec:revisedquestions}

If the knots $K_0$ and $K_1$ were $R$-concordant via an annulus $A$ in a CAT $4$-manifold $V$, then the meridians $\mu_0$ and $\mu_1$ would be freely homotopic in $V-A$. Hence, upon doing zero-framed surgery on $A$, one would have that $M_{K_0}$ is CAT $R$-homology cobordant to $M_{K_1}$ via a $4$-manifold $W$ with the additional property that the meridional elements are homologous in $W$, that is $(i_0)_*(\mu_0)=(i_1)_*(\mu_1)$ in $H_1(W;\Z)$. The rational homology cobordisms constructed in the proofs of Proposition~\ref{prop:mainresult} \textit{fail to have this additional property}. This suggests that the following revised question might be a better analogy to Question~\ref{question:KirbyQ}.

\begin{question}\label{revisedquestion} If $M_{K_0}$ and $M_{K_1}$ are CAT $\Q$-homology cobordant via a cobordism wherein the meridians are homologous (or freely homotopic), then are $K_0$ and $K_1$ CAT $\Q$-concordant?
\end{question}

Of course the answer is yes if one of the knots is the trivial knot.

One can revise the question in a different way. Following Kawauchi ~\cite{Ka4}, one can define $K_0$ and $K_1$ to be \textit{strongly $\Q$-concordant} if they are  $\Q$-concordant via an annulus $A$ in a $4$-manifold $V$ such that the complexity of the concordance is $1$ (see Section~\ref{sec:Qinvts}). This requires not only that the meridians are homologous, but also that they generate $H_1(V-A)$ modulo torsion.  Correspondingly one can say that $M_{K_0}$ and $M_{K_1}$ are \textit{strongly} CAT $\Q$-homology cobordant if they are CAT $\Q$-homology cobordant via $4$-manifold $W$ such that the merdians are are homologous and generate $H_1(W;\Z)$ modulo torsion.

\begin{question}\label{revisedquestion2} If $M_{K_0}$ and $M_{K_1}$ are strongly CAT $\Q$-homology cobordant, then are $K_0$ and $K_1$ strongly CAT $\Q$-concordant?
\end{question}

By using the techniques of Proposition~\ref{prop:Rbasichomcobordism} one can show that Question~\ref{revisedquestion2} also has a positive answer if one of the knots is the trivial knot.

Both of these questions are open as of this time.

\bibliographystyle{plain}
\bibliography{QConcBib}

\end{document}